\documentclass[12pt]{article}

\usepackage{amsmath,amsthm}

\theoremstyle{plain}
\newtheorem{lemma}{Lemma}[section]
\newtheorem{thm}[lemma]{Theorem}
\newtheorem{cor}[lemma]{Corollary}

\theoremstyle{definition} 
\newtheorem*{define}{Definition}

\newcommand{\ex}{\mathrm{ex}}
\newcommand{\gex}{\mathrm{gex}}
\newcommand{\N}{\mathbf{N}}

\title{Extensions of the linear bound in the F\"uredi--Hajnal conjecture}
\author{
Martin Klazar\thanks{Institute for Theoretical Computer Science (ITI) and Department of Applied Mathematics (KAM), 
Faculty of Mathematics and Physics of Charles University, Malostransk\'e n\'am\v est\'\i\ 25, 118 00 Praha, 
Czech Republic. Email: klazar@kam.mff.cuni.cz. 
ITI is supported by the project 1M0021620808 of the Czech Ministry of Education.
}
\and 
Adam Marcus\thanks{Research performed as Visiting Researcher, Alfr\'ed R\'enyi Institute, 1364 Budapest, Pf.127,
Hungary, on leave from the Department of Mathematics (ACO), Georgia Institute of Technology,
Atlanta, GA 30332-0160. Email: adam@math.gatech.edu.  
This research was made possible due to funding by the Fulbright Program in Hungary.}
}
\date{\today}

\begin{document}

\maketitle

\begin{abstract}
We present two extensions of the linear bound, due to Marcus and Tardos, on the number of 1's in an $n \times n$ 0-1 matrix avoiding a fixed permutation matrix. 
We first extend the linear bound to hypergraphs with ordered vertex sets and, using previous results of Klazar, we prove an exponential bound on the number of hypergraphs on $n$ vertices which avoid a fixed permutation.
This, in turn, solves various conjectures of Klazar as well as a conjecture of Br\"and\'en and Mansour.
We then extend the original F\"uredi--Hajnal problem from ordinary matrices to $d$-dimensional matrices and show that the number of 1's in a $d$-dimensional 0-1 matrix with side length $n$ which avoids a $d$-dimensional permutation matrix is $O(n^{d-1})$. 
  
\end{abstract}

\section{Introduction}\label{intro}

F\"uredi and Hajnal asked in \cite{fure_hajn} whether for every fixed 0-1 permutation matrix $P$ the maximum number 
of 1's in an $n\times n$ 0-1 matrix $M$ avoiding $P$ is $O(n)$; the avoidance here means that $P$ cannot 
be obtained from $M$ by a series of row deletions, column deletions, and replacements of $1$'s with $0$'s 
(in particular, permuting rows or columns of $M$ is not allowed).
The F\"uredi--Hajnal conjecture was settled by Marcus and Tardos in \cite{marc_tard} where they proved 
that if $M$ avoids a $k \times k$ permutation matrix, then the number of 1's in $M$ is at most $2k^4\binom{k^2}{k}n$.
In this paper we present extensions of their linear bound to more general structures.

The Marcus--Tardos bound can be reformulated in the language of graph theory, since matrices with entries 0 and 1 can be viewed as the incidence matrices of bipartite graphs. 
Thus if $P=([2k],E(P))$ is a graph on the vertex set $[2k]=\{1,2,\dots,2k\}$ with $k$ mutually disjoint edges, each of which connects the sets $[k]$ and 
$[k+1,2k]=\{k+1,k+2,\dots,2k\}$, and $M=([2n],E(M))$ is a graph on $[2n]$ which only has edges connecting 
$[n]$ and $[n+1,2n]$ and does not contain $P$ as an ordered subgraph, then $M$ has only linearly 
many edges, i.e. $|E(M)|=O(n)$. 
It is easy to modify the proof in \cite{marc_tard} so that it gives
a linear bound for all $P$-avoiding graphs $G$ (not necessarily bipartite) on the vertex set $[2n]$ (and therefore $[n]$). 
In Section 2 we extend this bound further to hypergraphs with edges of arbitrary size. We also discuss 
exponential enumerative bounds which follow from the linear extremal bounds as corollaries.

In yet another light, 0-1 matrices can be viewed as the (characteristic matrices of) binary relations. 
In Section 3 we generalize the original proof of Marcus and Tardos to $d$ dimensions and show that every $d$-ary relation on $[n]$ which avoids a fixed 
$d$-dimensional permutation has at most $O(n^{d-1})$ elements.

\section{Extensions to hypergraphs}

For a graph $G'=([k],E')$, we define $\gex_<(n,G')$ to be the maximum number $|E|$ of edges in a 
graph $G=([n],E)$ that does not contain $G'$ as an ordered subgraph. 
We represent a permutation $\pi=a_1a_2\dots a_k$ of $[k]$ by the graph
$$
P(\pi)=([2k],\{\{i,k+a_i\}:\ i\in[k]\}). 
$$
As we mentioned in Section~\ref{intro}, it is easy to modify the proof in \cite{marc_tard} to obtain the bound
\begin{equation}
\label{graph_extr_bound}
\gex_<(n,P(\pi))=O(n)
\end{equation}
where the constant in $O$ depends only on $\pi$.

For the hypergraph extension we need a few more definitions. 
A {\em hypergraph} is a finite collection $H=(E_i:\ i\in I)$ of finite nonempty edges 
$E_i$ which are subsets of $\N=\{1,2,\dots\}$. The {\em vertex set} is $V(H)=\bigcup_{i\in I}E_i$. For 
simplicity we do not allow (unlike in the graph case) isolated vertices; for our extremal problems 
this restriction is immaterial (isolated vertices in graphs can be represented by singleton edges in our extension).
In general we will allow multiple edges, and will denote a hypergraph as {\em simple} if it has no multiple edges. 
We say that $H'=(E_i':\ i\in I')$ is {\em contained} in $H=(E_i:\ i\in I)$, 
written $H'\prec H$, 
if there exists an increasing injection $f: V(H')\to V(H)$ and an injection $g:I'\to I$ such that 
$f(E_i')\subset E_{g(i)}$ for every $i\in I'$; otherwise we say that $H$ {\em avoids} $H'$. 
To put it differently, $H'\prec H$ means that $H'$ can be obtained from $H$ by deleting some edges,
deleting vertices from the remaining edges, and relabeling the vertices so that their ordering is 
preserved. This containment generalizes the ordered subgraph relation. Note that a simple hypergraph can contain a non-simple hypergraph. 

The {\em order} $v(H)$ of $H$ is the number of vertices $v(H)=|V(H)|$, the {\em size} $e(H)$ is 
the number of edges $e(H)=|I|$, and the {\em weight} $i(H)$ is the number of incidences 
$i(H)=\sum_{i\in I}|E_i|$. We define two hypergraph extremal functions. 

\begin{define}
Let $F$ be any hypergraph. We associate with $F$ the functions 
$\ex_e(\cdot,F),\ex_i(\cdot,F): \N\to\N$,
\begin{eqnarray*}
\ex_e(n,F)&=&\max\{e(H):\  H\not\succ F\;\&\;\mbox{$H$ is simple}\;\&\;v(H)\le n\}\\
\ex_i(n,F)&=&\max\{i(H):\  H\not\succ F\;\&\;\mbox{$H$ is simple}\;\&\;v(H)\le n\}.
\end{eqnarray*}
\end{define}
Obviously, $\ex_e(n,F)\le \ex_i(n,F)$ for every $F$ and $n$. 
If $F$ has at least two edges and has no two separated edges (edges $E_1$ and $E_2$ satisfying $E_1<E_2$), Klazar's Theorem~2.3 in \cite{klaz04} gives an inequality in the opposite direction: 
$$
\ex_i(n,F)\le (2v(F)-1)(e(F)-1)\ex_e(n,F).
$$
So in particular, for every permutation $\pi$ of $[k]$, 
\begin{equation}
\label{exi_from_exe}
\ex_i(n,P(\pi))\le(4k-1)(k-1)\ex_e(n,P(\pi)). 
\end{equation}
Thus a linear bound on $\ex_i(n,P(\pi))$ follows directly from one on $\ex_e(n,P(\pi))$. 

The latter bound can be derived using the techniques in \cite{klaz04} along with the graph bound in (\ref{graph_extr_bound}). 
To explain the reduction we need the notion of the {\em blow-up} of a graph. 
A graph $G'$ is an $m$-blow-up of a graph $G$ if for every edge coloring 
of $G'$ by colors from $\N$ such that every color is used at most $m$ times, there exists a subgraph of $G'$ that is order-isomorphic to $G$ and no two of its edges have the same color. 
Let $G$ be a graph with $k$ vertices, $H$ be a $\binom{k}{2}$-blow-up of $G$, and $f:\N\to\N$ be a function such that $\gex_<(n,H)<nf(n)$ for every $n\in\N$. 
Then Theorem~3.1 in \cite{klaz04} states that, for every $n\in\N$, 
\begin{equation}
\label{hgfromg}
\ex_e(n,G)<e(G)\cdot \gex_<(n,G)\cdot\ex_e(2f(n)+1,G).
\end{equation}
Combining the bounds in (\ref{graph_extr_bound}), (\ref{exi_from_exe}), and (\ref{hgfromg}) we obtain the following result:

\begin{thm}\label{C4}
For every permutation $\pi$,
$$
\ex_i(n,P(\pi))=O(n).
$$
\end{thm}
\begin{proof}
For $m\in\N$ and a $k$-permutation $\pi$, consider a permutation graph $P(\pi')$ that arises
from $P(\pi)$ by replacing every edge in $P(\pi)$ with a bundle of $k(m-1)+1$ edges so that the initial vertices of the 
edges in each bundle form an interval in $P(\pi')$ and the same holds for the final vertices. The positions of the intervals 
are as in $P(\pi)$, that is, for every selection of one edge from each bundle the resulting graph is 
order-isomorphic to $P(\pi)$. There are many such graphs $P(\pi')$  
($\pi'$ is a $(k^2(m-1)+k)$-permutation) and each of them is, by the pigeonhole principle, an $m$-blow-up of $P(\pi)$. 

We set $m=\binom{2k}{2}$. 
By the graph bound in (\ref{graph_extr_bound}), there are constants 
$c_{\pi}$ and $c_{\pi'}$ such that
$$
\gex_<(n,P(\pi))<c_{\pi}n\ \mbox{ and }\gex(n,P(\pi'))<c_{\pi'}n
$$
for every $n$. 
We set $H=P(\pi')$ and $f(n)=c_{\pi'}$ and apply the bound in (\ref{hgfromg}) to get the linear bound
$$
\ex_e(n,P(\pi))<kc_{\pi}\cdot\ex_e(2c_{\pi'}+1,P(\pi))\cdot n.
$$
By the bound in (\ref{exi_from_exe}), 
$$
\ex_i(n,P(\pi))<k(k-1)(4k-1)c_{\pi}\cdot\ex_e(2c_{\pi'}+1,P(\pi))\cdot n,
$$
proving the claim.
\end{proof}

Klazar posed the following six extremal and enumerative conjectures in \cite{klazar00}:

\begin{enumerate}
\item[C1:] The number of simple $H$ such that $v(H) = n$ and $H \not\succ P(\pi)$ is $< c_1^n$.
\item[C2:] The number of maximal simple $H$ with $v(H) = n$ and $H \not\succ P(\pi)$ is $< c_2^n$.
\item[C3:] For every simple $H$ with $v(H)=n$ and $H \not\succ P(\pi)$, we have $e(H) < c_3n$.
\item[C4:] For every simple $H$ with $v(H)=n$ and $H \not\succ P(\pi)$, we have $i(H) < c_4n$.
\item[C5:] The number of simple $H$ with $i(H) = n$ and $H \not\succ P(\pi)$ is $< c_5^n$.
\item[C6:] The number of $H$ with $i(H) = n$ and $H \not\succ P(\pi)$ is $< c_6^n$.
\end{enumerate}

He showed that all six conjectures hold for a large class of permutations $\pi$ and that they hold for every $\pi$ in weaker forms: with almost 
linear and almost exponential bounds (respectively). 
Conjecture~C4, however, is precisely the statement of Theorem~\ref{C4}, and it is easy to extend the proof given in this paper to affirm that all six conjectures hold for every permutation $\pi$.

We shall show how to amend the proofs in \cite{klazar00} to prove C1, and then note that C1 implies C2, C3, C5 and C6 via Lemma~2.1 of \cite{klazar00}.
\begin{cor}\label{C5C6}
For every permutation $\pi$ there exists a constant $c_1 > 1$ (depending on $\pi$) so that the number of simple hypergraphs on the vertex set $[n]$ avoiding $P(\pi)$ is $< c_1^n$.
\end{cor}
\begin{proof}
Theorems~2.4 and 2.5 in \cite{klazar00} show that the number of hypergraphs with a given weight $i(H)$ that avoid $P(\pi)$ is at most $9^{(3^{2k}+2k)i(H)}$.
Thus by Theorem~\ref{C4}, we are done.
\end{proof}

The Stanley--Wilf conjecture (see B\'ona \cite{bona 2004}), proved by Marcus and Tardos in \cite{marc_tard} as a corollary 
of their linear extremal bound, claimed that for every permutation $\pi$ there is a constant $c=c(\pi)$ such that  the number of permutations $\sigma$ 
of $[n]$ avoiding $\pi$ is $<c^n$; the avoidance of permutations here means that $P(\sigma)$ is not an ordered subgraph of $P(\pi)$. 
In view of the reformulation from permutations to bipartite graphs mentioned in Section~\ref{intro}, Corollary~\ref{C5C6} is an extension of the Stanley--Wilf conjecture. 
A related extension was proposed by Br\"and\'en and Mansour in Section~5 of \cite{bran_mans}: they conjectured that the number of 
{\em words} over the ordered alphabet $[n]$ which have length $n$ and avoid $\pi$ is at most exponential in $n$. 
These words can be represented by simple graphs $G$ on $[2n]$ in which every edge connects $[n]$ and $[n+1,2n]$ and 
every $x \in [n]$ has degree exactly $1$; the containment of ordered words is then just the ordered subgraph relation. 
Hence this extension is subsumed in Corollary~\ref{C5C6}. 

Corollary~\ref{C5C6} subsumes yet another extension of the Stanley--Wilf conjecture to set partitions proposed 
by Klazar \cite{klaz_eujc00}. 
This extension is related to $k$-{\em noncrossing} and $k$-{\em nonnesting} set partitions whose 
exact enumeration was recently investigated by Chen et al. \cite{chen_al} and Bousquet-M\'elou and Xin \cite{bm_xin}. 
Consider, for a set partition $H$ of $[n]$, the graph $G(H)=([n],E)$ in which an edge connects two neighboring elements of 
a block (not separated by another element of the same block). Thus $H$ is represented by increasing paths which 
are spanned by the blocks. $H$ is a $k$-noncrossing (resp. $k$-nonnesting) partition iff $P(12\dots k)$ 
(resp. $P(k(k-1)\dots 1)$) is not an ordered subgraph of $G(H)$.
Thus Corollary~\ref{C5C6} provides an exponential bound: for fixed $k$, the numbers of $k$-noncrossing and $k$-nonnesting 
partitions of $[n]$ grow at most exponentially. 

\section{An extension to $d$-dimensional matrices}\label{adam}

We now generalize the original F\"uredi--Hajnal conjecture from ordinary 0-1 matrices to $d$-dimensional 0-1 matrices. 
As was mentioned in Section~\ref{intro}, these are just $d$-ary relations (or, as we will discuss later, $d$-uniform, 
$d$-partite hypergraphs).
We keep the matrix terminology, however, both for the sake of consistency and to highlight the similarities with the 
original Marcus--Tardos proof in \cite{marc_tard}.

\begin{define}
We will call a $(d+1)$-tuple $M = (M; n_1,\dots,n_d)$ where $M \subset [n_1] \times \ldots \times [n_d]$ a 
$d$-{\em dimensional (0-1) matrix}, and will refer to the elements of $M$ as {\em edges}. 
\end{define}

If $F=(F;k_1,\dots,k_d)$ and $M=(M;n_1,\dots,n_d)$ are two $d$-dimensional matrices, we say that $F$ is {\em contained} in 
$M$, $F \prec M$, if there exist $d$ increasing injections $f_i:[k_i] \to [n_i]$, $i=1,2,\dots,d$, such that for 
every $(x_1, \dots, x_d) \in F$ we have $(f_1(x_1), \dots, f_d(x_d))\in M$; otherwise we say that $M$ {\em avoids} $F$.

\begin{define}
We set $f(n,F,d)$ to be the maximum size $|M|$ of a $d$-dimensional matrix 
$(M;n,\dots,n)$ that avoids a $d$-dimensional matrix $F$. 
\end{define}

For $i \in [d]$, we will denote the {\em projection mapping} from $[n_1] \times\dots \times [n_d]$ to $[n_i]$ as $\pi_i$. 
For $t \in [d]$, we define the $t$-{\em remainder} of $M=(M;n_1,\dots,n_d)$ to be the $(d-1)$-dimensional matrix 
$N=(N;n_1',\dots,n_{d-1}')$ where $n_1'=n_1,\dots,n_{t-1}'=n_{t-1}$, $n_t'=n_{t+1},\dots,n_{d-1}'=n_d$ and the edge 
$(e_1,\dots,e_{d-1})\in N$ if and only if 
$(e_1,\dots,e_{t-1},x,e_t,e_{t+1},\dots,e_{d-1})\in M$ for some $x\in[n_t]$. 

Let $I_1 < I_2 < \dots < I_r$ be a partition of $[n]$ into $r$ intervals and $M=(M;n,\dots,n)$ a $d$-dimensional matrix. 
We define the {\em contraction} of $M$ (with respect to the intervals) to be the $d$-dimensional matrix $N=(N;r,\dots,r)$ 
given by 
$(e_1,\dots, e_d)\in N$ iff $M\cap (I_{e_1} \times \dots \times I_{e_d})\ne\emptyset$ (we could define the contraction 
operation 
for a general $d$-dimensional matrix and with distinct and general partitions in each coordinate but we will 
not need such generality).

We say that $P=(P;k,\dots,k)$ is a $d$-{\em dimensional permutation of} $[k]$
if for every $i \in [d]$ and $x \in [k]$ there is exactly one edge $e \in P$ with $\pi_i(e)=x$. 
Note that $|P|=k$ and that there are exactly $(k!)^{d-1}$ $d$-dimensional permutations of $[k]$. 
For $d=1$, the only 1-dimensional permutation $(P;k)$ is $[k]$, and for $d=2$ the 2-dimensional permutations 
$P=(P;k,k)$ are precisely the $k \times k$ 0-1 permutation matrices. 
A $d$-dimensional permutation of $[k]$ can be 
thought of as a $d\times k$ matrix with the first row normalized to $1,2,\dots,k$ and with each row being a permutation of 
$1,2,\dots,k$. The columns would then give the coordinates of the $k$ edges in $P$.

It is also possible to view the structure $M=(M;n_1, \ldots, n_d)$ as an ordered, $d$-uniform, 
$d$-partite hypergraph with partitions $[n_i]$.  In this interpretation, the image of $M$ by the projection $\pi_i$ 
is obtained by intersecting the edges with the $i^{th}$ partition, while the intersections with the union of all partitions 
except the $t^{th}$ one give the $t$-remainder of $M$ (in both cases we disregard multiplicity of edges).
Furthermore, the set of $d$-dimensional permutations of $[k]$ would be the set of perfect matchings of the complete $d$-uniform, 
$d$-partite hypergraph on $kd$ vertices.

We will make use of two observations, analogous to those made in \cite{marc_tard}:
\begin{enumerate}
\item For any $t \in [d]$, the $t$-remainder of a $d$-dimensional permutation of $[k]$ is a $(d-1)$-dimensional permutation of $[k]$.
Furthermore, each edge of the resulting $t$-remainder can be completed (by adding the $t$-th coordinate) in a unique way to an edge of 
the original permutation. 
\item If $M=(M;n,\dots,n)$ avoids a $d$-dimensional permutation, then so does any contraction of $M$.
\end{enumerate}

\begin{thm}
\label{thm:dd}
For every fixed $d$-dimensional permutation $P$, $$f(n,P,d) = O(n^{d-1}).$$
\end{thm}

\noindent
On the other hand it is clear that for a $d$-dimensional permutation $P$ with $|P| > 1$ we have $f(n,P,d)\geq n^{d-1}$
($f(n,P,d)=0$ if $|P|=1$). Thus, for a $d$-dimensional permutation $P$ with $|P|>1$, 
$$
f(n,P,d)=\Theta(n^{d-1}).
$$
This bound can be given an equivalent formulation. We say that a matrix $M=(M;n_1,\dots,n_d)$ is a $d$-{\em dimensional}
$k$-{\em grid} if each $[n_i]$ can be partitioned in $k$ intervals 
$I_{i,1}<I_{i,2}<\dots<I_{i,k}$ so that $|M\cap(I_{1,j_1}\times I_{2,j_2}\times\cdots\times I_{k,j_k})|=1$ for every 
$d$-tuple $(j_1,j_2,\dots,j_d)\in[k]^d$ (thus, in particular, $|M|=k^d$). Let $g(n,k,d)$ be the maximum size of a 
$d$-dimensional $n\times n\times\dots\times n$ matrix that contains no $d$-dimensional $k$-grid. Then
$$
g(n,k,d)=\Theta(n^{d-1}).
$$
It is clear that $g(n,k,d)\ge n^{d-1}$. The bound $g(n,k,d)=O(n^{d-1})$ implies $f(n,P,d)=O(n^{d-1})$ for every $P$ because
every $d$-dimensional $k$-grid contains every $d$-dimensional permutation of $[k]$. In the other way, it is easy to see 
that there exist $d$-dimensional $k$-grids that are $d$-dimensional permutations of $[k^d]$. Thus $f(n,P,d)=O(n^{d-1})$ 
implies $g(n,k,d)=O(n^{d-1})$. 

To prove Theorem~\ref{thm:dd}, we will show that a $d$-dimensional matrix of big enough size must contain every $d$-dimensional 
permutation of $k$.
We set 
$$f(n,k,d)=\max_P f(n,P,d)$$ 
where $P$ runs through all $d$-dimensional permutations of $[k]$.

\begin{lemma}
\label{lem:ineq}
Let $d \geq 2$, $m,n_0 \in \N$.  Then 
$$
f(mn_0,k,d)\leq (k-1)^d\cdot f(n_0,k,d) + dn_0m^{d} \binom{m}{k}\cdot f(n_0,k,d-1).$$
\end{lemma}
\begin{proof}
Let $M=(M,mn_0, \ldots, mn_0)$ be a $d$-dimensional matrix that avoids $P$, a $d$-dimensional permutation of $[k]$.  
We aim to bound the size of $M$.

We split $[mn_0]$ into $n_0$ intervals $I_1<I_2<\dots<I_{n_0}$, each of length $m$, and define, 
for $i_1,\dots,i_d\in[n_0]$, 
$$
S(i_1,\dots,i_d)=\{e \in M : \pi_j(e) \in I_{i_j}\mbox{ for }j=1,\dots,d\}.
$$
Note that this partitions the set of edges of $M$ into $n_0^d$ pieces.
We will call these sets of edges {\em blocks} and we define a cover of these blocks by a total of $dn_0 + 1$ sets 
$\{U_0\} \cup \{U(t,j) : t\in[d],j\in[n_0]\}$ as follows:
\begin{itemize}
\item $U(t,j)= \{S(i_1,\dots,i_d):\ i_t=j \mbox{ and } |\pi_t(S(i_1,\dots,i_d))|\geq k\}$
\item $U_0$ consists of the blocks which are not in any $U(t,j)$
\end{itemize}

Note that the total number of non-empty blocks is exactly the number of edges in the contraction of $M$ with respect to the partition 
$\{I_i\}$.  
Since $M$ does not contain $P$, the contraction of $M$ can not contain $P$, so the number of non-empty blocks is at most $f(n_0, k, d)$.  
Also note that any block $B$ in $U_0$ has at most $(k-1)^d$ edges in it (because $B\subset X_1\times\dots\times X_d$ for some 
$X_i\subset[mn_0]$ with $|X_i|<k$). Hence 
$${\textstyle
|\bigcup U_0|\le (k-1)^d \cdot f(n_0, k, d).
}
$$ 

Now we fix $t\in [d]$ and $j \in [n_0]$. Clearly,
$${\textstyle
|\bigcup U(t,j)|\le m^d|U(t,j)|.
}
$$
We assume, for a contradiction, that $|U(t,j)| > \binom{m}{k} \cdot f(n_0,k,d-1)$. 
By the definition of $U(t,j)$ and the pigeonhole principle, there are $k$ numbers $c_1<c_2<\dots<c_k$ in $I_j$ and $r$ blocks 
$S_1,S_2,\dots,S_r$ in $U(t,j)$ where $r > f(n_0,k,d-1)$ such that for every $S_a$ and every $c_b$ there is an $e\in S_a$ with 
$\pi_t(e)=c_b$. 
Let $P'$ be the $t$-remainder of $P$ and $M'=(M';n_0,\dots,n_0)$ be the $(d-1)$-dimensional matrix arising from contracting 
$(\bigcup_{i=1}^r S_i,n,\dots,n)$ with respect to the intervals $I_i$ and then taking the $t$-remainder. 
Since $|M'| = r > f(n_0,k,d-1)$, $M'$ contains $P'$. 
Thus among the blocks $S_1,S_2,\dots,S_r$ there exist $k$ of them---call them $S_1,S_2,\dots,S_k$---so that for any selection of $k$ edges 
$e_1 \in S_1,\dots,e_k\in S_k$ their $t$-remainders form a copy of $P'$. 
Furthermore, due to the property of the blocks $S_i$, it is possible to select $e_1,\dots,e_k$ so that their $t$-th coordinates agree 
with $P$. Then $e_1,\dots,e_k$ form a copy of $P$, a contradiction. Therefore
$${\textstyle
|\bigcup U(t,j)| \leq m^d|U(t,j)|\le m^d\binom{m}{k} \cdot f(n_0,k,d-1)
}
$$
and
$${\textstyle
|\bigcup_{t,j} U(t,j)|\le dn_0m^d\binom{m}{k} \cdot f(n_0,k,d-1).
}
$$
Combining this with the bound for $U_0$ gives the stated bound.
\end{proof}

Theorem~\ref{thm:dd} will be a direct consequence of the following lemma:
\begin{lemma}
\label{lem:induct} For $m = \lceil k^{d/(d-1)} \rceil$, $f(n,k,d) \leq k^d\Big(dm\binom{m+1}{k}\Big)^{d-1}n^{d-1}$.
\end{lemma}
\begin{proof}
We will proceed by induction on $d + n$.
For $d=1$ this holds since $f(n,k,1) = k-1$ and for $n < k^2$, this holds trivially. 
Now given $n$ and $d \geq 2$, assume that the hypothesis is true for all $d',n'$ such that $d' + n' < d + n$. 

Let $n_0=\lfloor n/m \rfloor$ and 
$$c_d = k^d\bigg(dm\binom{m+1}{k}\bigg)^{d-1}.$$
Using the inequality $f(n,k,d) < f(mn_0,k,d) + dmn^{d-1}$, Lemma~\ref{lem:ineq}, the inductive hypotheses, and 
$n_0 \leq n/m$, we get
$$f(n,k,d)<\Bigg(\frac{(k-1)^d}{m^{d-1}}c_d+dm\bigg(\binom{m}{k}c_{d-1}+1\bigg)\Bigg)n^{d-1}.$$
Since $\frac{(k-1)^d}{m^{d-1}}\le(1-\frac{1}{k})^d \leq 1-\frac{1}{k}$ and $\binom{m}{k}c_{d-1}+1 
\leq \binom{m+1}{k}c_{d-1}$, it follows that $f(n,k,d)<c_dn^{d-1}$ with the above defined $c_d$.
\end{proof}

\section{Concluding remarks}
We were informed recently that Balogh, Bollob\'as and Morris \cite{balo_boll_morr} derived Theorem~\ref{C4} (their Theorem~2) and Corollary~\ref{C5C6} (their Theorem~1) independently.
The proofs in \cite{balo_boll_morr} are self-contained (not appealing to the results in \cite{klaz04}) and their approach is different from ours.
In fact, they are able to prove stronger statements, which in turn imply Theorem~\ref{C4} and Corollary~\ref{C5C6} from this paper.

It should be noted that we make no effort to optimize any of the constants in Section~\ref{adam}, however it would be interesting to see if any of the constants could be drastically reduced.  
The constant achieved in this paper is double exponential in $k$, whereas we conjecture the true constant is in fact much smaller.

\section{Acknowledgments}
The authors would like to thank G\'abor Tardos for, among other things, his enlightening conversations and helpful remarks.

\end{document}